\documentstyle[12pt]{article}
\newcommand{\mysection}[1]{\section{#1}\setcounter{equation}{0}}

\begin{document}%

\LARGE
\begin{center}
{\bf The asymptotic distribution of the largest prime divisor}
\end{center}

\normalsize


\newtheorem{sub}{\name}[section]
\newtheorem{subn}{\name}
\newtheorem{Thm}{Theorem}[section]
\newtheorem{Lem}[Thm]{Lemma}
\newtheorem{Prop}[Thm]{Proposition}
\newtheorem{Cor}[Thm]{Corollary}
\newtheorem{Rem}[Thm]{Remark}
\newtheorem{Def}[Thm]{Definition}
\newtheorem{Ex}[Thm]{Example}
\renewcommand{\thesubn}{}
\newcommand{\dn}[1]{\def\name{#1}}
\newcommand{\be}{\begin{equation}}
\newcommand{\ee}{\end{equation}}
\newcommand{\bs}{\begin{sub}}
\newcommand{\es}{\end{sub}}
\newcommand{\bsn}{\begin{subn}}
\newcommand{\esn}{\end{subn}}
\newcommand{\bea}{\begin{eqnarray}}
\newcommand{\eea}{\end{eqnarray}}
\newcommand{\BA}[1]{\begin{array}{#1}}
\newcommand{\EA}{\end{array}}

\newcommand{\Real}{\mbox{${\rm I\!R}$}}
\newcommand{\real}{{\rm I\!R}}
\newcommand{\Nat}{\mbox{${\rm I\!N}$}}
\newcommand{\thkl}{\rule[-.5mm]{.3mm}{3mm}}
\newcommand{\Proof}{\mbox{\noindent {\bf Proof} \hspace{2mm}}}
\newcommand{\mbinom}[2]{\left (\!\!{\renewcommand{\arraystretch}{0.5}
 \mbox{$\begin{array}[c]{c}
 #1\\ #2
 \end{array}$}}\!\! \right )}
\newcommand{\brang}[1]{\langle #1 \rangle}
\newcommand{\vstrut}[1]{\rule{0mm}{#1mm}}
\newcommand{\rec}[1]{\frac{1}{#1}}
\newcommand{\set}[1]{\{#1\}}
\newcommand{\dist}[2]{\mbox{\rm dist}\,(#1,#2)}
\newcommand{\opname}[1]{\mbox{\rm #1}\,}
\newcommand{\supp}{\opname{supp}}
\newcommand{\mb}[1]{\;\mbox{ #1 }\;}
\newcommand{\undersym}[2]
 {{\renewcommand{\arraystretch}{0.5}
 \mbox{$\begin{array}[t]{c}
 #1\\ #2
 \end{array}$}}}

\newlength{\wex}  \newlength{\hex}
\newcommand{\understack}[3]{%
 \settowidth{\wex}{\mbox{$#3$}} \settoheight{\hex}{\mbox{$#1$}}
 \hspace{\wex}
 \raisebox{-1.2\hex}{\makebox[-\wex][c]{$#2$}}/
 \makebox[\wex][c]{$#1$}
  }%

\newcommand{\smit}[1]{\mbox{\small \it #1}}
\newcommand{\lgit}[1]{\mbox{\large \it #1}}
\newcommand{\scts}[1]{\scriptstyle #1}
\newcommand{\scss}[1]{\scriptscriptstyle #1}
\newcommand{\txts}[1]{\textstyle #1}
\newcommand{\dsps}[1]{\displaystyle #1}

\def\ga{\alpha}     \def\gb{\beta}       \def\gg{\gamma}
\def\gc{\chi}       \def\gd{\delta}      \def\ge{\epsilon}
\def\gth{\theta}                         \def\vge{\varepsilon}
\def\gf{\phi}       \def\vgf{\varphi}    \def\gh{\eta}
\def\gi{\iota}      \def\gk{\kappa}      \def\gl{\lambda}
\def\gm{\mu}        \def\gn{\nu}         \def\gp{\pi}
\def\vgp{\varpi}    \def\gr{\rho}        \def\vgr{\varrho}
\def\gs{\sigma}     \def\vgs{\varsigma}  \def\gt{\tau}
\def\gu{\upsilon}   \def\gv{\vartheta}   \def\gw{\omega}
\def\gx{\xi}        \def\gy{\psi}        \def\gz{\zeta}
\def\Gg{\Gamma}     \def\Gd{\Delta}      \def\Gf{\Phi}
\def\Gth{\Theta}
\def\Gl{\Lambda}    \def\Gs{\Sigma}      \def\Gp{\Pi}
\def\Gw{\Omega}     \def\Gx{\Xi}         \def\Gy{\Psi}
\mysection{Introduction}\label{section1}
A point ${\bf z}$ in $R^m$ is {\em a lattice point} if ${\bf z}=(z_1,...,z_m)$ where each
$z_j$ is an integer. Consider the number of
lattice points included in the simplex $S(a_1,...,a_m)$, where
\be \label{simplex}
S(a_1,...,a_m)=\left \{{\bf z}:\sum _{j=1}^m\frac {z_j}{a_j}\leq 1,\,z_j\geq 0,
1\leq j\leq m \right \},
\ee
and $a_j$, $j=1,2,...,m$, are positive real numbers. Denote this number by
$\rho (a_1,...,a_m)$, or $\rho (S)$.

We need estimates of $\rho (S)$ as a tool in studying the following problem.
Let $n$ and $N$ be two positive real numbers, and we are interested in the number of integers
$2\leq k\leq N$ such that the largest prime divisor of $k$ does not exceed $n$. We
denote this number by
\be \label{Nzero}
\nu(n,N).
\ee
Denote by $\{p_j\}_{j=1}^{\infty }$ the
increasing sequence of the primes, and let $m$ be such that
\be \label{mdefn}
p_{m}<n\leq p_{m+1}.
\ee
Then by the Prime Numbers Theorem
\be \label{maprox}
m\approx \frac {n}{\ln n}
\ee
in the sense that the ratio between the two sides of (\ref{maprox}) tends to 1 as
$n\to\infty $. We are thus interested in the integers $k\leq N$ which are of the form
\be \label{kform}
k=\prod _{j=1}^mp_j^{x_j},\,x_j \mbox { are nonnegative integers}.
\ee
Equivalently, we are interested in integers $k$ as in
(\ref{kform}) for which
\be \label{simformula}
\sum _{j=1}^m(\ln p_j)x_j\leq \ln N
\ee
holds, and we have to estimate
\be \label{lnpjsimp}
\rho \left (\frac {\ln N}{\ln p_1},...,\frac{\ln N}{\ln p_m}\right ).
\ee

Concerning $\nu (n,N)$ we have the following result, which is a corollary of our
main results, Theorems \ref{mainlwrbndrslt8} and \ref{mainuprbndrslt8}. It deals
with situations where
\be \label{intrnge}
\ln N<<n<<N,
\ee
in a sense expressed precisely in the theorem.
\begin{Thm} \label{mainres0}
$(i)$ Consider pairs $(n,N)$ such that
$$
\frac {\ln \ln N}{\ln n}\to 0 \mbox { and }
\frac {\ln N}{\ln n}\to \infty \mbox { as } n\to \infty .
$$
Then
\be \label{lwrboundassy11}
\frac {\ln \nu (n,N)}{\ln N}>1-\frac {\ln \ln N+\ln \ln \ln N}
{\ln n}+\frac {a^{\star }
+\ln \ln n+\delta (n,N)}{\ln n},
\ee
where $a^{\star }=1+\ln (e-1)$ and $\delta (n,N)\to \infty $ as $n\to \infty $.
\newline
$(ii)$
Consider pairs $(n,N)$ such that
$$
N<e^{\sqrt n} \mbox { and }\frac {\ln N}{(\ln n)^2\ln \ln N}\to \infty
\mbox { as } n\to \infty ,
$$
and let $a$ be any constant. Then there exists $n_0$ such that
\be \label{nineqult31}
\frac {\ln \nu (n,N)}{\ln N}<1-\frac {\ln \ln N+\ln \ln \ln N}{\ln n}
+\frac {a+\ln \ln n+\ln \ln \ln n}{\ln n}
\ee
if $n>n_0$.
\end{Thm}
\begin{Rem}\label{gapconj}
There is a gap between the lower bound $(\ref{lwrboundassy11})$ and the upper bound
$(\ref{nineqult31})$, where in the former we have $a^{\star }+\delta (n,N)$ while in the latter
$a+\ln \ln \ln n$. The following might consist of sharper bounds. For an integer $k\geq 2$
and sufficiently large $N$ denote
$$
\ln ^{(k)}N=\ln \cdots \ln N
$$
where the logarithm function appears $k$ times, and denote
$$
\mbox {\rm Ln}^{(k)}N=\sum _{j=2}^k\ln ^{(j)}N.
$$
We conjecture that the sharper bounds
$$
\frac {\ln \nu (n,N)}{\ln N}>1-\frac {\mbox {\rm Ln}^{(k)}N}
{\ln n}+\frac {a^{\star }+\mbox {\rm Ln}^{(k-1)}n+\delta (n,N)}{\ln n},
$$
and
$$
\frac {\ln \nu (n,N)}{\ln N}<1-\frac {\mbox {\rm Ln}^{(k)}N}{\ln n}
+\frac {a+\mbox {\rm Ln}^{(k)}n}{\ln n}
$$
may be established. We note that these bounds reduce to $(\ref{lwrboundassy11})$ and
$(\ref{nineqult31})$ for $k=3$.
\end{Rem}
The following result covers a different range of parameters $N$ and $n$.
\begin{Thm}\label{niceresult}
Consider the set $E$ of integers $1\leq k\leq N$ for which all the
prime divisors are smaller than $\sqrt N$. In our notations $\#(E)=\nu (\sqrt N,N)$, and
we have that
\be \label{FsqrtN}
\nu (\sqrt N,N)>\alpha N \mbox { for some constant }\alpha >0 \mbox { and every }N>1.
\ee
Actually, for sufficiently large $N$ we may take $\alpha =\ln (e/2)$ in $(\ref{FsqrtN})$.
\end{Thm}
The proof is relegated to the appendix.

The next result will be needed below.
\begin{Lem} \label{thelemma}
The following relation holds:
\be \label{estimate1}
\rho (a_1,...,a_m)> \frac {\prod _{j=1}^ma_j}{m!}.
\ee
\end{Lem}
{\em Proof}: The proof is by induction on $m$. For $m=1$ we have
$$\rho (a_1)=[a_1]+1>a_1,$$
so that (\ref{estimate1}) holds in this case. (We denote by $[x]$ the integer part of
$x$.)

Let $m\geq 2$ and assume that the assertion of the proposition holds for $m-1$. Denote
$$a_m=a$$
and
$$
\rho(a_1,...,a_{m-1})=\rho _0.
$$
Let $0\leq j\leq [a]$ be an integer, and we consider the $(m-1)$-simplex
$$
S_j=S(a_1,...,a_m)\cap \{x_m=j\}.
$$
Then $S_j$ is a translation of the $(m-1)$-simplex
$$
S[(1-j/a)a_1,...,(1-j/a)a_{m-1}],
$$
and by the induction hypothesis, the number of lattice points in $S_j$, denoted $\rho _j$,
satisfies
$$
\rho _j>\left (1-\frac{j}{a}\right )^{m-1}\frac {\prod _{j=1}^{m-1}a_j}{(m-1)!}.
$$
Since
$$
\rho (a_1,...,a_m)=\sum _{j=0}^{[a]}\rho _j,
$$
it follows that
\be \label{kappasum}
\rho (a_1,...,a_m)>\frac {\prod _{j=1}^{m-1}a_j}{(m-1)!}
\sum _{j=0}^{[a]}\left (1-\frac{j}{a}\right )^{m-1},
\ee
and we estimate the sum in (\ref{kappasum}) by an integral as follows:
\be \label{ntesti}
\sum _{j=0}^{[a]}\left (1-\frac{j}{a}\right )^{m-1}>
\int _0^{a}\left (1-\frac{x}{a}\right )^{m-1}dx=\frac {a}{m}.
\ee
Using (\ref{ntesti}) in (\ref{kappasum}) implies (\ref{estimate1}), concluding
the proof.$\hfill \Box$

For parameters in a certain range the estimate of $\rho $ in
(\ref{estimate1}) is adequate,
while for others it is quite poor. For example, consider the situation where
$a_j=L$ for every $1\leq j\leq m$, in which case (\ref{estimate1}) yields the lower
bound $L^m/m!$. Assuming that $m$ is large, we use Stirling's formula
\be \label{strlngformla}
m!\approx \sqrt{2\pi m}\left (\frac {m}{e}\right )^m
\ee
 to approximate
\be \label{Stirling}
\frac {L^m}{m!}\approx \frac {1}{\sqrt{2\pi m}}\left (\frac {eL}{m}\right )^m.
\ee
This yields a good estimate if $1<<m<L$, but it provides a very poor bound if, e.g.,
$L<m/2e$. In this case (\ref{Stirling}) yields a lower bound which is smaller
than $2^{-m}$, while actually $\rho $ increases to infinity as $L\to \infty $.

Considering (\ref{lnpjsimp}) we take as a typical order of magnitude for $a_j$ in
(\ref{estimate1})
$$
\frac {\ln N}{\ln p_j}\approx \frac {\ln N}{\ln n},
$$
and the above discussion implies that the estimate in (\ref{estimate1}) is
poor if
$$\frac {\ln N}{\ln n}<\frac {m}{2e},$$
namely
\be \label{badestim}
n>2e\ln N.
\ee
In the range of parameters (\ref{intrnge}) which is of interest for us
(\ref{badestim}) certainly holds, and the estimate in (\ref{estimate1})
is actually useless.
\begin{Rem}\label{nopolynom}
If $n=(\ln N)^q$ then by Theorem $\ref{mainres0}$, $\nu (n,N)\approx N^{1-1/q}$.
Therefore, if $n$ is of a polynomial order in $\ln N$,
then the set of integers having largest prime factor that is smaller than $n$ is sparse
in $[2,N]$. On the other hand, if $n> \sqrt N$ then by Theorem $\ref{niceresult}$
$$\nu (n,N)>\alpha N$$
for some constant $\alpha >0$, so that
the set of integers having largest prime factor in $[2,n]$ is quite dense in $[2,N]$.
It is thus of interest to study the situation where $\ln \ln N<<\ln n<<\ln N$.
\end{Rem}

The paper is organized as follows. In the next section we establish a preliminary
bound, which will be improved in the sequel. In section \ref{section3}
we describe a setting which enables  the study of tight lower and upper bounds for
$\nu(n,N)$. In section \ref{section4} we introduce a family
of auxiliary problems in which our problem can be imbedded. In section \ref{section5}
we introduce our iterations method, which is the main technical tool developed in
this paper. In sections \ref{section6} and \ref{section7} we establish lower and upper
bounds for the auxiliary problems, and our main results are presented in section
\ref{section8}. In the appendix we establish Theorem \ref{niceresult} and
Proposition \ref{upperbundln2}.
\mysection{A Preliminary lower bound for $\nu(n,N)$}\label{section2}

To compute a lower bound for $\nu (n,N)$ we will estimate the number of lattice points
which are contained in the simplex (\ref{simformula}) (where $m$ is as in (\ref{mdefn})
and (\ref{maprox})). Since $\ln p_j<\ln n$, it follows that this number is larger
than the number of lattice points contained in the simplex
\be \label{nsimplex}
\sum _{j=1}^mx_j\leq \frac {\ln N}{\ln n},\,x_j\geq 0.
\ee
Obviously, the number $\displaystyle{\sum _{j=1}^mx_j}$
is an integer whenever $(x_1,...,x_m)$
is a lattice point. Hence the number of lattice points contained in the simplex
(\ref{nsimplex}) is equal to
\be \label{vsum}
\sum _{k=1}^{l}f(k,m),
\ee
where
\be \label{nunot}
l=\left [\frac {\ln N}{\ln n}\right ],
\ee
and where $f(k,m)$ denotes the number of different ways in which $k$ can be
written as a sum of $m$ nonnegative integers. Clearly
\be \label{fkdefn}
f(k,m)=\left (\begin{array}{cc}k+m-1\\k\end{array}\right )=
\frac {m(m+1)\cdots (m+k-1)}{k!},
\ee
so that the number of lattice points contained in the
simplex (\ref{simformula}) is larger than
\be \label{Sksum}
\sum _{k=1}^{l}\left (\begin{array}{cc}k+m-1\\k\end{array}\right ).
\ee
We express the $k$th term in (\ref{Sksum}) in the form
\be \label{fkmexprss}
\left (\begin{array}{cc}k+m-1\\k\end{array}\right )=\frac {m^k}{k!}
\left (1+\frac {1}{m}\right )\cdots \left (1+\frac {k-1}{m}\right ),
\ee
and it follows that
$$
f(k,m)>\frac {m^k}{k!}
$$
for every $k$. By (\ref{vsum}), the quantity $m^k/k!$
is a lower bound for $\nu (n,N)$ for each $1\leq k\leq l$ , and
we note that if $m>>l$ (namely $n>>\ln N$), then the lower bound $m^l/l!$
is much larger than the lower bound $l^m/m!$ which results
from (\ref{estimate1}).

Using $f(l,m)$ as a lower bound for $\nu $ and employing
Stirling's formula (\ref{strlngformla}) we  obtain
\be \label{N0est}
\nu(n,N)> \frac {1}{\sqrt {2\pi l}}\left (\frac {em}{l}\right )^{l}
>\frac {1}{\sqrt {\ln N}} \left (\frac {en}{\ln N}\right )^{\ln N /\ln n}
\ee
if $n>n_0$ for some $n_0$.
In case that $l$ is large enough so that Stirling's approximation
$(\ref{strlngformla})$ may be employed for it, then (\ref{N0est}) may be expressed
in the form
\be \label{lowrbound}
\frac {\ln \nu(n,N)}{\ln N}>1-\frac{\ln \ln N}{\ln n}+\frac {1}{\ln n}
-\frac {\ln \ln N}{2\ln N}.
\ee

To obtain upper bounds for $\nu (n,N)$ the following result will be useful.
\begin{Prop}\label{nurecrse}
Let $\{p_k\}_{k=1}^{\infty }$ denote the sequence of primes. Then
\be \label{nupkrecurs}
\nu (p_{k+1},N)=\sum _{j=0}^{\left [\frac {\ln N}{ln p_{k+1}}\right ]}\nu
\left (p_k,N/p_{k+1}^j\right )
\ee
holds for every $N>2$ and $k\geq 1$.
\end{Prop}
{\em Proof}: Let ${\cal F}_k(N)$ denote the set of integers $z\leq N$ whose largest
prime divisor does not exceed $p_k$, so that
\be \label{nuFkNrl}
\nu (p_k,N)=\#\{{\cal F}_k(N)\}.
\ee
Denote by $A_j$ the set of integers $z\in {\cal F}_{k+1}(N)$ such that $p_{k+1}^j$
is the largest power of $p_{k+1}$ which divides $z$.
It is then easy to see that
\be \label{Ajdefinit}
A_j=p_{k+1}^j{\cal F}_k\left (\frac {N}{p_{k+1}^j}\right )
\ee
and
\be \label{Fk1union}
{\cal F}_{k+1}(N)=\bigcup _{j\geq 0}A_j,
\ee
a disjoint union. The relation (\ref{nupkrecurs}) follows
from (\ref{nuFkNrl}), (\ref{Ajdefinit}) and (\ref{Fk1union}). $\hfill \Box$

We obtain the following result, which will be used in section \ref{section7}.
\begin{Prop}\label{upperbundln2}
Let $\alpha >0$ be fixed, and consider pairs $(n,N)$ such that
\be \label{nalpN2}
n=\alpha (\ln N)^2.
\ee
Then there exists a constant $C>1$ such that
\be \label{nlnNsqrbnd}
\frac {\ln \nu (n,N)}{\ln N}< 1-\frac {\ln \ln N}{\ln n}+\frac {C}{\ln n}
\ee
holds for every $N>1$, where $n$ is as in $(\ref{nalpN2})$.
\end{Prop}
The proof is displayed in the appendix.

\mysection{The reduced order simplex}\label{section3}
In this section we relate with the high dimensional simplex (\ref{simformula}) a
simplex of smaller order. We will study certain properties of this simplex,
which will be used in the next sections as tools used to establish tight lower
and upper bounds for the number of solutions of (\ref{simformula}).

In establishing a lower bound in section \ref{section2} we used the inequality
\be \label{simplnpjn}
\ln p_j<\ln n
\ee
for every $1\leq j\leq m$. Modifying this approach
we divide the integers interval $(1,n)$ into subintervals
\be \label{Jidfn}
J_i=\left (\frac {n}{e^i},\frac {n}{e^{i-1}}\right ), i=1,2,...,r,
\ee
where
\be \label{r1}
r= [\ln n] \mbox { if } \ln n<[\ln n]+\ln 2
\ee
and
\be \label{r2}
r= [\ln n]+1 \mbox { if } \ln n>[\ln n]+\ln 2.
\ee
For simplicity of notations we henceforth consider only case (\ref{r1}),
and comment that the discussion and main results in case (\ref{r2}) are the same.
(In Remark \ref{rdiffrnc} we will indicate where the difference between case (\ref{r1})
and case (\ref{r2}) plays a role.)

Refining (\ref{simplnpjn}) we have for primes $p_j\in J_i$ the relations
\be \label{r3}
\ln n-i<\ln p_j<\ln n-i+1,
\ee
and regarding (\ref{simformula}) this implies
\be \label{Jisum}
(\ln n-i)z_i<\sum _{p_j\in J_i}(\ln p_j)x_j<(\ln n-i+1)z_i,
\ee
where we denote
\be \label{Aidefn}
z_i=\sum _{p_j\in J_i}x_j.
\ee
Clearly $(z_1,...,z_r)$ is a nonnegative lattice point in $R^r$.
\begin{Rem}\label{rdiffrnc}
The cases $(\ref{r1})$ and $(\ref{r2})$ differ only when considering $i=r$
in the left hand side of $(\ref{r3})$.
\end{Rem}

If $\{x_j\}_{j=1}^{m}$ is a solution of (\ref{simformula}), then in
view of (\ref{Jisum}) this implies
\be \label{Aiconsrnt}
\sum _{i=1}^r(\ln n-i)z_i<\ln N.
\ee
Therefore the number of solutions
$\{x_j\}_{j=1}^m$ of (\ref{simformula}) is smaller than the number of solutions
$\{x_j\}_{j=1}^m$ of (\ref{Aiconsrnt}).
(We say that $\{x_j\}_{j=1}^m$ is a solution of (\ref{Aiconsrnt}) if
(\ref{Aidefn}) and $(\ref{Aiconsrnt})$ are satisfied.)
Similarly, if $\{x_j\}_{j=1}^m$ is a solution of
\be \label{Aircnt}
\sum _{i=1}^r(\ln n-i+1)z_i<\ln N,
\ee
then in view of (\ref{Jisum}) it is also a solution of (\ref{simformula}),
implying that the number of solutions
$\{x_j\}_{j=1}^m$ of (\ref{simformula}) is larger than the number of solutions
$\{x_j\}_{j=1}^m$ of (\ref{Aircnt}). These considerations are the basis of our
computation of upper and lower bounds for $\nu(n,N)$.

For a prescribed lattice point $(z_1,...,z_r)$ which satisfies (\ref{Aiconsrnt}) we
are interested in the number of lattice points $\{x_j\}_{j=1}^m$ in $R^m$ for which
(\ref{Aidefn}) holds for every $i=1,2,...,r$. Let $m_i$ denote the size of the set
$\{j:p_j\in J_i\}$:
$$
m_i=\#\left \{p_j\in \left (\frac {n}{e^{i}},\frac {n}{e^{i-1}}\right )\right \},
$$
and if $m_i>>1$, then by the Prime Numbers Theorem
\be \label{miapproxi}
m_i\approx \frac {(e-1)n}{(\ln n-i)e^{i}},
\ee
and we have the inequality
\be \label{midefn}
m_i> \frac {n}{e^{i}(\ln n-i)}.
\ee
Employing the notation $f(k,m)$ in (\ref{fkdefn}), it follows that
the number of lattice points $\{x_j\}_{j=1}^m$ that satisfy (\ref{Aidefn})
for every $1\leq i\leq r$ is
\be \label{Kexpress}
K(z_1,...z_r)=\prod _{i=1}^rf(z_i,m_i).
\ee

We Denote by $\overline {\nu }(n,N)$ and $\underline {\nu }(n,N)$ the
number of solutions of (\ref{Aiconsrnt}) and (\ref{Aircnt}) respectively,
and it follows that $\nu (n,N)$ is bounded from above
by  $\overline {\nu }(n,N)$ and from below by $\underline {\nu }(n,N)$.
Using the expression $K(z_1,...,z_r)$ in (\ref{Kexpress}) we consider sums of the form
\be \label{Mexpress}
M(F)=\sum _{{\bf z}\in F}K(z_1,...,z_r),
\ee
where the summation runs over all the lattice points ${\bf z}=\{z_1,...,z_r\}$ which belong to
some set $F$ in $R^r$. Thus when $F$ in (\ref{Mexpress}) is the set of points belonging to
the simplex (\ref{Aiconsrnt}), denoted $F_1$, then by (\ref{fkdefn}) and (\ref{Kexpress})
we have
\be \label{upperbnd}
\overline {\nu}(n,N)=\sum _{\{z_i\}\in F_1}\prod _{i=1}^{r}\frac {m_i^{z_i}}{z_i!}
\left (1+\frac {1}{m_i}\right )\cdots \left (1+\frac {z_i-1}{m_i}\right ).
\ee
Similarly we obtain the following lower bound for $\nu$
\be \label{lowerbnd}
\underline {\nu }(n,N)=\sum _{\{z_i\}\in F_2}\prod _{i=1}^{r}\frac {m_i^{z_i}}{z_i!},
\ee
where $F_2$ is the set of all the lattice points in the simplex (\ref{Aircnt}).

We next consider the product
$$
P_i=\prod _{k=1}^{z_i-1}\left (1+\frac {k}{m_i}\right )
$$
that appears in the right hand side of (\ref{upperbnd}),
and in view of the inequality $\ln (1+x)<x$ for $x>0$ we obtain
$\ln P_i<z_i^2/2m_i$,
hence
\be \label{Pibound}
P_i<e^{z_i^2/2m_i}.
\ee
When dealing with a lower bound we will ignore
the term $\prod _{i=1}^{r}P_i$ in the right hand side of (\ref{upperbnd}), and
we will focus on computing a lower bound to expressions of the form
\be \label{Zexprs}
Z(F)=\sum _{\{z_i\}\in F}\prod _{i=1}^{r}\frac {m_i^{z_i}}{z_i!}
\ee
for certain sets $F$. We will then describe the modifications required to obtain
an upper bound by taking into consideration the terms $P_i$ in (\ref{upperbnd}).

\mysection{A family of auxiliary problems}\label{section4}
It will be convenient to study our main problem, of estimating sums of the form
(\ref{Mexpress}), by using slightly different notations.
In this section we define a collection of problems, parameterized by two real variables,
such that for certain values of the parameters the auxiliary problem coincides with the
main problem. Thus for a positive number $c>1$, let $r=[c]$ and consider the inequality
\be \label{cconstrnt}
cz_0+(c-1)z_1+(c-2)z_2+\cdots +(c-r+1)z_{r-1}<M
\ee
for some positive number $M>1$, where ${\bf z}=\{z_i\}_{i=0}^{r-1}$ is a nonnegative
lattice point in  $R^{r}$ (compare with (\ref{Aircnt})).
We associate with $c$ the $r$ bases
\be \label{mbases}
m_i=\frac {(e-1)e^{c-i}}{c-i},\,0\leq i\leq r-1
\ee
(compare with  (\ref{miapproxi}) in case that $c=\ln n$). In view of (\ref{lowerbnd})
we address the problem of computing the sum
\be \label{cMsum}
F(c,M)=\sum _{\bf z}\prod _{i=0}^{r-1}\frac {m_i^{z_i}}{z_i!},
\ee
where ${\bf z}=(z_0,...,z_{r-1})$ runs over all the nonnegative lattice points which satisfy
(\ref{cconstrnt}); we call this {\em Problem $P_{c,M}$} for the
$r$ variables $z_0$,...,$z_{r-1}$.
\begin{Rem}\label{ourapplic}
There is a close relation between the value of Problem $P_{c,M}$ and $\nu (n,N)$ for
\be \label{ourvalus}
c=\ln n\mbox { and } M=\ln N.
\ee
Thus the value of $P_{c,M}$ yields a lower bound for $\nu(n,N)$. We also note that
if $c\geq M$ $($namely $n\geq N)$ and $N$ is an integer, then
\be \label{nuNMeql}
\nu (n,N)=N=e^M.
\ee
\end{Rem}

To establish an upper bound for $\nu (n,N)$ we will estimate a sum of the type
(\ref{Mexpress}), which is associated with the simplex
\be \label{inter}
(c-1)z_1+(c-2)z_2+\cdots +(c-r)z_{r}<M
\ee
(compare with (\ref{Aiconsrnt})). This sum is smaller than the corresponding sum that
is associated with the simplex
\be \label{ccconstrnt}
cz_0+(c-1)z_1+(c-2)z_2+\cdots +(c-r)z_{r}<M,
\ee
which we denote by $G_0(c,M)$. Thus to obtain an upper bound for $G_0(c,M)$
we consider a sum similar to the one in (\ref{cMsum}), where
we take into consideration the
terms $P_i$ in (\ref{Pibound}). We then address the problem of computing the sum
\be \label{GcMsum}
G(c,M)=\sum _{\bf z}\prod _{i=0}^{r}\frac {m_i^{z_i}e^{z_i^2/m_i}}{z_i!},
\ee
where ${\bf z}=(z_0,z_1,...,z_{r})$ runs over all the nonnegative lattice points which
satisfy (\ref{ccconstrnt}); we call this {\em Problem $Q_{c,M}$} for the
$r+1$ variables $z_0$,$z_1$,...,$z_{r}$.
\begin{Rem}
We use the simplex $(\ref{ccconstrnt})$ rather than the simplex $(\ref{inter})$, which is
more directly related to $(\ref{Aiconsrnt})$, in order to avoid repetition of
computations for the lower and upper bounds. Thus a substantial part of the computations
for $(\ref{cconstrnt})$ and $(\ref{ccconstrnt})$ will be unified.
\end{Rem}

We claim that for a fixed value of $z_0$, Problem $P_{c,M}$ reduces to
Problem $P_{c-1,M-cz_0}$ for the $r-1$ variables $z_1$,...,$z_{r-1}$. To justify
this statement we have to check that the $r-1$ bases $m_1$,...,$m_{r-1}$ in
(\ref{mbases}) are indeed the bases associated with Problem $P_{c-1,M-cz_0}$,
which is easily verified.

The possible values for the variable $z_0$ in (\ref{cconstrnt}) are the
integers $z$ satisfying
$$
0\leq z\leq \frac {M}{c},
$$
and it follows from (\ref{cMsum}) that
\be \label{Frecurs}
F(c,M)=\sum _{z=0}^{[M/c]}F(c-1,M-cz)\frac {m_0^z}{z!}.
\ee

In the subsequent discussion we will consider situations where $F(\cdot ,\cdot )$
satisfies inequalities of the form
\be \label{Fform}
F(c,M)\geq Be^{M\left (1-\frac {\ln M}{c+1}
+\frac {\gamma }{c+1}\right )}
\ee
for some constant  $0<B\leq 1$. In terms of the original parameters we are actually
interested in inequalities of the form
\be \label{nuform}
\nu (n,N)\geq BN^{\left (1-\frac {\ln \ln N}{\ln n+1}
+\frac {\gamma }{\ln n+1}\right )},
\ee
where $(n,N)$ and $(c,M)$ are related as in (\ref{ourvalus}).
\begin{Rem} \label{Mc2}
It follows from $(\ref{FsqrtN})$ in Theorem $\ref{niceresult}$
that for a fixed $\gamma $, inequality $(\ref{nuform})$
holds whenever $M/c<2$. Indeed, for $M=\ln N$ and
$c=\ln n$ the condition $M/c<2$ translates to $n>\sqrt N$, and
$\nu (n,N)>\alpha N$ by $(\ref{FsqrtN})$. But the inequality
$$
\alpha N>N^{1-\frac {\ln \ln N}{\ln n+1}+\frac {\gamma }{\ln n+1}}
$$
is equivalent to
$$
\frac {\ln N}{\ln n+1}(\ln \ln N-\gamma )>-\ln \alpha ,
$$
and this holds for every $N>N_0$, for some $N_0$, since $n<N$. For $N\leq N_0$,
however, $(\ref{nuform})$ holds for some $B(\gamma )$,
since in this case we have a bounded set of pairs $(n,N)$.
Therefore, when trying to establish an inequality of the type $(\ref{Fform})$,
we may assume that
\be \label{Mclrg2}
\frac {M}{c}\geq 2,
\ee
since for $M/c<2$ inequality $(\ref{nuform})$ is already established.
\end{Rem}
\mysection{The iterations method}\label{section5}
The discussion in this section is fundamental to our analysis. We develop the iterations
method which will be employed in the subsequent sections to establish lower and upper
bounds for $\nu $.

Assume that for a certain $\gamma >0$ and some $0<B<1$,
inequality (\ref{Fform}) holds for any pair $(c,M)$ which verifies
\be \label{cMM0cond}
c\leq \kappa _0
\ee
for a certain $\kappa _0$. We consider then pairs $(c,M)$ that satisfy
\be \label{nextblock1}
\kappa _0<c\leq \kappa _0+1,
\ee
and our goal is to establish the inequality $(\ref{Fform})$ for such pairs as well. Once
this is achieved we will iterate the argument to obtain a lower bound for all pairs in a
certain domain.

Intending to employ (\ref{Frecurs}) to establish a lower bound to $F(c,M)$, and assuming
that (\ref{Fform}) holds whenever (\ref{cMM0cond}) is satisfied, we
will estimate from below the expressions
\be \label{Fzbarexp}
F(c-1,M-cz)\frac {m_0^{z}}{z!}
\ee
for integers $0\leq z\leq M/c$. By (\ref{nextblock1}) $c-1\leq \kappa _0$, and
we may use (\ref{Fform}) for the pair $(c-1,M-cz)$, obtaining
\be \label{c1Mzbnd}
F(c-1,M-cz)\geq Be^A,
\ee
where
\be \label{Alogaritm}
A=(M-cz)-\frac {1}{c}(M-cz)\ln (M-cz)+\frac {(M-cz)\gamma}{c}.
\ee
Also
\be \label{mzeE}
\frac {m_0^z}{z!}>e^E,
\ee
denoting
\be \label{Elogaritm1}
E=\left (z\ln m_0-z\ln z+z\right )-\left (\frac {1}{2}\ln z
+\frac {1}{2}\ln \pi +\frac {3}{2}\ln 2\right ),
\ee
where we used Stirling's formula
\be \label{Stirapprx}
St(z)=\sqrt {2\pi z}\left (\frac {z}{e}\right )^z
\ee
to estimate
\be \label{2gamaz}
z!<2St(z) \mbox { for every } z\geq 1.
\ee
A term $(-\ln 2)$ in (\ref{Elogaritm1}) arises from the factor
$2$ in (\ref{2gamaz}), and the term
\be \label{dist}
-\frac {1}{2}(\ln z+\ln \pi +\ln 2)
\ee
in (\ref{Elogaritm1}) is due to the logarithm of $\sqrt {2\pi z}$ in (\ref{Stirapprx}).
To avoid the disturbing term (\ref{dist}) in (\ref{Elogaritm1}) we note that
\be \label{zpibeta}
\frac {1}{2}\ln z+\frac {1}{2}\ln \pi +\frac {3}{2}\ln 2<\beta z
\ee
where $\beta >0$ may be chosen arbitrarily
small provided that $z$ is sufficiently large. It follows that
\be \label{zineqmodi}
z-\left (\frac {1}{2}\ln z
+\frac {1}{2}\ln \pi +\frac {3}{2}\ln 2\right )>bz
\ee
where
\be \label{b1beta}
b=1-\beta
\ee
may be chosen arbitrarily close to 1
provided that $z$ is large enough, and we thus obtain
\be \label{Elogaritm}
E>\left (z\ln m_0-z\ln z+bz\right )
\ee
for sufficiently large values of $z$.

It follows from $m_0=(e-1)e^{c}/c$ that
$$
z\ln m_0=cz-z\ln c+z\ln (e-1).
$$
Using the last equation in (\ref{Elogaritm}) and recalling (\ref{Alogaritm}) yield that
\be \label{AEestim}
A+E>H(z),
\ee
denoting
\be \label{maxexpres}
H(z)=M\left (1+\frac {\gamma }{c}\right )+(a-\gamma )z-\frac {M}{c}\ln c
-z\ln z-\left (\frac {M}{c}-z\right )\ln \left (\frac {M}{c}-z\right )
\ee
and
\be \label{adefinitn}
a=b+\ln (e-1).
\ee
Thus $a$ is smaller and arbitrarily close to $a^{\star }$, which is defined by
\be \label{limita}
a^{\star }=1+\ln (e-1).
\ee
It follows from (\ref{c1Mzbnd}), (\ref{mzeE}) and (\ref{AEestim}) that
\be \label{FcMHD}
F(c-1,M-cz)\frac {m_0^z}{z!}>Be^{H(z)},
\ee
and to obtain a lower bound for the sum in (\ref{Frecurs}) we will estimate
the maximal value of $H(z)$, $0\leq z\leq [M/c]$, where $z$ is an integer.

\begin{Rem}\label{newrem}
We will compute a maximizer $z_0$ of $H(\cdot )$ defined on the real interval $[0,[M/c]]$,
and in general $z_0$ is not an integer. Let $z_1$ be the integer
$$
z_1=z_0+\theta \mbox { for some }
0\leq \theta <1,
$$
and then
$$
H(z_1)=H(z_0)+\frac {1}{2}H''(\zeta )\theta ^2
$$
for some $z_0<\zeta <z_1$. But
$$
H''(\zeta )=\frac {-M/c}{\zeta (M/c-\zeta )},
$$
and it follows from $\zeta \geq 1$ that
$$
|H''(\zeta )|\leq \frac {M/c}{M/c-1}<2
$$
$($since $M/c>2)$, and we obtain
\be \label{Hz0z1}
H(z_1)>H(z_0)-\theta ^2.
\ee
Similarly, for the integer $z_2=z_0-(1-\theta )$ we have
\be \label{Hz0z2}
H(z_2)>H(z_0)-(1-\theta )^2.
\ee
It follows from $(\ref{FcMHD})$, $(\ref{Hz0z1})$ and $(\ref{Hz0z2})$ that
\be \label{intgrestm}
\sum _{z=0}^{[M/c]}F(c-1,M-cz)\frac {m_0^z}{z!}>B\left (e^{H(z_1)}+e^{H(z_2)}\right )
>Be^{H(z_0)}
\ee
since
$$
\min _{0\leq \theta \leq 1}\left \{e^{-\theta ^2}+e^{-(1-\theta )^2}\right \}>1.
$$
Therefore we may use the maximal value of $H(z)$ over the whole real interval
$0\leq z\leq M/c$.
\end{Rem}
We have the following basic result.
\begin{Prop}\label{maxHvalue}
Let $H(z)$ be as in $(\ref{maxexpres})$. Then
\be \label{Hzmaxlower1}
\max \left \{H(z):0\leq z\leq \frac {M}{c}\right \}=
M\left (1-\frac {\ln M}{c}+\frac {\gamma +f(\gamma )}{c}\right ),
\ee
where
\be \label{fdefn}
f(\gamma )=\ln (1+e^{a-\gamma }).
\ee
\end{Prop}
{\em Proof}:
Denoting
$$
K=\frac {M}{c} \mbox { and }z=Kt
$$
it follows that
\begin{eqnarray}\label{zmaxval}
\max _z\{(a-\gamma )z-z\ln z-(K-z)\ln (K-z)\}=\nonumber \\
-K\ln K+K\max _{0\leq t\leq 1}\{(a-\gamma )t-t\ln t-(1-t)\ln (1-t)\}.
\end{eqnarray}
We denote
\be \label{varphidfn}
\varphi (t)=(a-\gamma )t-t\ln t-(1-t)\ln (1-t),
\ee
and it follows that the maximizer $t_0$ of $\varphi $ satisfies
$$
(a-\gamma )-\ln t_0+\ln (1-t_0)=0.
$$
We conclude that
\be \label{toexpr}
t_0(\gamma )=\frac {1}{1+e^{\gamma -a}},
\ee
and the maximal value of $\varphi (\cdot )$ is given by
$$
(a-\gamma )t_0+\ln (1+e^{\gamma -a})-(1-t_0)(\gamma -a),
$$
which yields
\be \label{tmaxval}
\max \{\varphi (t):0\leq t\leq 1\}=\ln (1+e^{a-\gamma }).
\ee
We thus conclude from (\ref{zmaxval}) and (\ref{tmaxval}) that
\begin{eqnarray} \label{Hzmaximal}
\max _{0\leq z\leq K}\left \{z(a-\gamma )
-z\ln z-\left (\frac {M}{c}-z\right )
\ln \left (\frac {M}{c}-z\right )\right \}
=\nonumber \\ -K\ln K
+K\ln (1+e^{a-\gamma }).
\end{eqnarray}
It follows from (\ref{maxexpres}) and (\ref{Hzmaximal}) that (\ref{Hzmaxlower1}) is
satisfied, where $f(\gamma )$ in (\ref{fdefn}) is the maximum
in (\ref{tmaxval}). The proof of the proposition is complete. $\hfill \Box$
\newline
It follows from (\ref{Frecurs}), (\ref{intgrestm}) and (\ref{Hzmaxlower1}) that
\be \label{Hzmaxlower}
F(c,M)\geq B\exp \left \{M\left (1-\frac {\ln M}{c}
+\frac {\gamma+f(\gamma )}{c}\right )\right \}.
\ee

For the induction argument we need that (\ref{Fform}) would hold for some initial value of
$c$, say for $c=\kappa $ for some $\kappa >1$. This is the content of the following result.
\begin{Prop}\label{Bkapagama}
For a prescribed $\gamma >0$ the inequality
\be \label{Fkapgamform}
F(\kappa ,M)\geq B(\kappa ,\gamma )e^{M\left (1-\frac {\ln M}{\kappa +1}
+\frac {\gamma }{\kappa +1}\right )}
\ee
holds for every $M\geq 0$, where
\be \label{Bkapgamval}
B(\kappa ,\gamma )=e^{-e^{\kappa +\gamma }}.
\ee
\end{Prop}
{\em Proof}: The maximal value of
$$
M\mapsto M\left (1-\frac {\ln M}{c+1}+\frac {\gamma }{c+1}\right )
$$
is $ \displaystyle {\frac {e^{\kappa +\gamma }}{\kappa +1}}$, and it is attained at
$M_0=e^{\kappa +\gamma }$. Since $B(\kappa ,\gamma )$ in (\ref{Bkapgamval}) satisfies
$$
B(\kappa ,\gamma )e^{\frac {e^{\kappa +\gamma }}{\kappa +1}}<1,
$$
and since $F(c,M)\geq 1$, inequality (\ref{Fkapgamform}) follows for every $M>1$.
$\hfill \Box$.

We note that if $B$ is equal to $B(\kappa ,\gamma )$
in (\ref{Bkapgamval}), then (\ref{Fform}) holds for any pair $(c,M)$ such that
$c\leq \kappa$.
\mysection{A lower bound for Problem $P_{c,M}$}\label{section6}
In this section we employ the results of the previous section to establish a lower bound
for Problem $P_{c,M}$. We will construct a sequence
\be \label{cjMjsequ}
\{(c_j,M_j)\}_{j=0}^l
\ee
for which
(\ref{Hzmaxlower}) will be employed successively. The coefficient $B$ will be chosen
such that
\be \label{mainesti}
F(c,M)\geq B\exp \left \{M\left (1-\frac {\ln M}{c+1}
+\frac {\gamma '}{c+1}\right )\right \}
\ee
will hold for the pair $(c_l,M_l)$ for a certain $\gamma '=\gamma _l$, and consequently,
employing (\ref{Hzmaxlower}), it will
hold for each $(c_j,M_j)$ with a certain $\gamma '=\gamma _j$, in
particular for $(c,M)=(c_0,M_0)$.

Recall that in deriving the estimate (\ref{Hzmaxlower}) we used a value
$$
z_0=Kt_0
$$
which is associated with a pair $(c_1,M_1)$ such that $c_1=c_0-1$, and
\be \label{M1expr}
M_1=M_0(1-t_0).
\ee
Although it does not correspond to an integer $z$, it may be used to obtain
a lower bound for $F(c,M)$, as explained in Remark \ref{newrem}.

Concerning (\ref{Hzmaxlower}), we wish to estimate its right hand side as follows:
\be \label{Mgamampr}
M\left (1-\frac {\ln M}{c}
+\frac {\gamma+f(\gamma )}{c}\right )>
M\left (1-\frac {\ln M}{c+1}
+\frac {\gamma '}{c+1}\right )
\ee
for a certain $\gamma '$. Clearly the inequality (\ref{Mgamampr}) is equivalent to
\be \label{Mgamampr1}
\frac {\gamma +f(\gamma )}{c}>\frac {\ln M}{c(c+1)}+\frac {\gamma '}{c+1}.
\ee
For any $\beta >0$ we denote
\be \label{Dbetadfn}
D_{\beta }=\{(c,M):1\leq M\leq e^{\beta c}\},
\ee
and for a fixed $0<\alpha <a$ we denote for every pair $(c,M)$
\be \label{gamcM}
\gamma _{c,M}=a-\alpha +\ln c-\ln \ln M.
\ee

For a pair $(c,M)$ we consider the maximization over $z$ of
\be \label{maxc1Fczm0}
F(c-1,M-cz)\frac {m_0^z}{z!}.
\ee
We assume validity of (\ref{mainesti}) with $c-1$ replacing $c$, taking for $(c-1,M')$
$$
\gamma '=\gamma _{c-1,M-cz},
$$
namely we assume that
\be \label{Fc1lowbnd}
F(c-1,M')\geq B \exp \left \{M'\left (1-\frac {\ln M'}{c}
+\frac {\gamma _{c-1,M'}}{c}\right )\right \}
\ee
for every $1\leq M'\leq M$. Using (\ref{gamcM}) in (\ref{Fc1lowbnd}) yields
$$
F(c-1,M')\geq B \exp \left \{M'\left (1-\frac {\ln M'}{c}
+\frac {a-\alpha +\ln (c-1)-\ln \ln M}{c}\right )\right \},
$$
which we write in the form
\be \label{Fc1M'gam0}
F(c-1,M')\geq B \exp \left \{M'\left (1-\frac {\ln M'}{c}
+\frac {\gamma _0}{c}\right )\right \}
\ee
for every $1\leq M'\leq M$, denoting
\be \label{gama0expr}
\gamma _0=a-\alpha +\ln (c-1)-\ln \ln M.
\ee
The fact that the parameter $\gamma _0$ in (\ref{Fc1M'gam0}) is the same for all $M'$
enables to employ the results of section \ref{section5}. Thus the maximal value of
(\ref{maxc1Fczm0}) exceeds the maximal value which is obtained when we replace $F(c-1,M-cz)$
by the right hand side of (\ref{Fc1lowbnd}), with $M'=M-cz$, namely the maximal value of
\be \label{boundsmaxi}
\exp \left \{(M-cz)\left [1-\frac {\ln (M-cz)}{c}
+\frac {\gamma _0}{c}\right ]\right \}\frac {m_0^z}{z!}
\ee
over $0\leq z\leq M/c$. This latter maximum is attained at
\be \label{M'Mdfn}
M'=M(1-t_0)
\ee
where
\be \label{t0fdormul}
t_0=\frac {1}{1+e^{\gamma _0-a}}.
\ee
\begin{Prop}\label{Dbetacontn}
Let $\alpha >0$ be fixed. Then there exists a constant $\beta _0$ such that
\be \label{betaclose}
(c,M)\in {\cal D}_{\beta }\Rightarrow (c-1,M')\in {\cal D}_{\beta }
\ee
for every $0<\beta <\beta _0$.
\end{Prop}
{\em Proof}: By (\ref{gama0expr})
$$
e^{\gamma _0-a}=e^{-\alpha }\frac {c-1}{\ln M},
$$
and using this in (\ref{t0fdormul}) yields
\be \label{t0estabv}
t_0>e^{\alpha /2}\frac {\ln M}{c}
\ee
if
$$
\frac {\ln M}{c}<\beta _0
$$
for some $\beta _0$ which is small enough, and if $c$ is large enough. It follows
from (\ref{M'Mdfn}) that
$$
\ln M'<\ln M-t_0
$$
which, in view of (\ref{t0estabv}), yields
\be \label{lnMprlnM}
\ln M'<\ln M\left (1-\frac {e^{\alpha /2}}{c}\right ),
\ee
implying
\be \label{lnM'c1alph}
\frac {\ln M'}{c-1}<\frac {\ln M}{c}\left (\frac {c-e^{\alpha /2}}{c-1}\right ).
\ee
Thus (\ref{betaclose}) follows from (\ref{lnM'c1alph}), since $\alpha >0$. $\hfill \Box$

We will next establish (\ref{mainesti}) with
\be \label{gamaprchoic}
\gamma '=\gamma _{c,M}
\ee
(recall (\ref{gamcM})), assuming the validity of (\ref{mainesti}) with $c$ being replaced
by $c-1$.
\begin{Prop}\label{iteratineq}
Let $z_0$ be the maximizer in the maximization over $z$ of $(\ref{boundsmaxi})$, and let
$a$ be associated with $z_0$ as in $(\ref{zpibeta})$, $(\ref{b1beta})$ and
$(\ref{adefinitn})$. Let $\gamma '$ be as in $(\ref{gamaprchoic})$ and $\gamma =\gamma _0$
$($recall $(\ref{gama0expr}))$. Then $(\ref{Mgamampr})$ holds.
\end{Prop}
{\em Proof}: We consider the expression
\be \label{fgamacM}
f(\gamma )=f(\gamma _0)=\ln \left (1+e^{\alpha }\frac {\ln M}{c-1}\right ).
\ee

For any $0<q<1$, which may be arbitrarily close to $1$, we have that
\be \label{fgamma0est}
\ln \left (1+e^{\alpha }\frac {\ln M}{c-1}\right )>qe^{\alpha }\frac {\ln M}{c-1}
\ee
if $(\ln M)/(c-1)$ is sufficiently small.
But $\alpha >0$ is fixed while $q$ is arbitrarily close to $1$, and it follows
from (\ref{fgamma0est}) that there exist $c_0$ and $\beta $ such that
\be \label{fgammaestimt}
f(\gamma )>\frac {\ln M}{c}
\ee
if $c>c_0$ and $(c,M)\in {\cal D}_{\beta }$.

For $\gamma=\gamma _0 $ and $\gamma '$ as in (\ref{gama0expr}) and (\ref{gamaprchoic})
the inequality
\be \label{gamgampr}
\frac {\gamma }{c}>\frac {\gamma '}{c+1}
\ee
is equivalent to
\be \label{equivform}
a-\alpha -\ln \ln M+(c+1)\ln (c-1)>c\ln c
\ee
But (\ref{equivform}) follows from
$$
\ln c<\ln (c-1)+\frac {1}{c-1}
$$
in view of $M<e^{c-1}$. The inequality (\ref{Mgamampr}) is a consequence of
(\ref{Mgamampr1}), (\ref{fgammaestimt}) and (\ref{gamgampr}). $\hfill \Box$

For a fixed $\beta >0$ we have relation (\ref{betaclose}), which enables
to use (\ref{Mgamampr}) iteratively.
It follows from (\ref{Hzmaxlower}), (\ref{Mgamampr}) and Proposition \ref{iteratineq}
that for a fixed $\alpha >0$, the inequality
\be \label{lowerboundcM}
F(c,M)>B\exp \left \{M\left (1-\frac {\ln M}{c+1}+\frac {a-\alpha +\ln c-\ln \ln M}{c+1}
\right )\right \}
\ee
holds for certain pairs $(c,M)$. More precisely, the above discussion yields
the next iterative property.
\begin{Prop}\label{kapkap1lwr}
For a fixed $\alpha >0$ there exist $\kappa _0>0$ and $\beta >0$ with the following property:
If $\kappa >\kappa _0$ is such that  $(\ref{lowerboundcM})$ holds for every
$(c,M)\in {\cal D}_{\beta }$ satisfying $\kappa _0<c\leq \kappa $, then it also
holds for every $(c,M)$ that verifies
$$
(c,M)\in {\cal D}_{\beta } \mbox { and } \kappa _0<c\leq \kappa +1.
$$
\end{Prop}
\begin{Rem}\label{aastar}
Consider a sequence $(\ref{cjMjsequ})$ where $(c_{j-1},M_{j-1})$ is the maximizing pair
associated with $(c_j,M_j)$ in the above discussion. We denote by $t_j$, $z_j$ and $a_j$ the
corresponding parameters in this maximization, and it follows from $(\ref{t0estabv})$ that
$$
t_j>\frac {\ln M_j}{c_j}.
$$
Then the maximizer $z_j$ satisfies
$$
z_j=\frac {M_j\ln M_j}{c_j},
$$
and in view of $(\ref{lnM'c1alph})$ it follows that $z_j\to \infty $
if $M_j\to \infty $. But then by $(\ref{zpibeta})$, $(\ref{b1beta})$ and
$(\ref{adefinitn})$,
we may take $a_j\to a^{\star }$, since $\alpha >0$ may be chosen arbitrarily
small. We conclude that if $M_j\ln M_j/c_j\to \infty $ for the sequence $(\ref{cjMjsequ})$
then we may assume that $a_j\to a^{\star }$.
\end{Rem}
To start the iterations procedure we need the following result:
\begin{Prop}\label{lwrboundres}
For a fixed $\alpha >0$ let $\kappa _0$ and $\beta $ be as in Proposition
$\ref{kapkap1lwr}$, and let $B$ be defined by
\be \label{Bkap0}
B=e^{-\kappa _0e^{a+\kappa _0}}.
\ee
Then $(\ref{lowerboundcM})$ holds for every $(c,M)\in {\cal D}_{\beta }$ such that
$c\geq \kappa _0$.
\end{Prop}
{\em Proof}: The assertion of the proposition follows from Propositions \ref{Bkapagama}
and \ref{kapkap1lwr}, employing an induction argument. $\hfill \Box$

We conclude from Propositions \ref{kapkap1lwr} and \ref{lwrboundres} the following result.
\begin{Prop}\label{noaymptot}
Let $a<a^{\star }$ be fixed. Then there exist $\beta >0$, $c_0$ and $B$ such that
\be \label{lwrboundassy}
F(c,M)>B\exp \left \{M\left (1-\frac {\ln M+\ln \ln M}{c+1}+\frac {a+\ln c}{c+1}
\right )\right \}
\ee
for every $(c,M)$ such that $M<e^{\beta c}$ and $c>c_0$.
\end{Prop}

The following is the asymptotic lower bound which we obtain for $F(c,M)$.
By Remark \ref{aastar} we may assume that $a$ is arbitrarily close to $a^{\star }$,
provided that $M(\ln M)/c$ is sufficiently large. We therefore may
replace $a$ and $\alpha $ in (\ref{lowerboundcM}) by $a^{\star }+\delta (c,M)$, where
$\delta \to 0$ if $M(\ln M)/c\to \infty $. Moreover, we note that the denominator
$c+1$ in (\ref{lwrboundassy}) may be replaced by $c$, as expressed in (\ref{lwrboundassy1}),
since the difference that arises from this change may be absorbed into a term $\delta (c,M)$
as in (\ref{lwrboundassy1}) and (\ref{deltacMlwr}). We further note that
the coefficient $B$ in (\ref{lwrboundassy}) may be absorbed
in $\delta (c,M)$ under the assumption $M/c\to \infty $.
\begin{Thm} \label{mainlwrbndrslt}
Consider pairs $(c,M)$ such that
\be \label{lwrcondi}
\frac {\ln M}{c}\to 0 \mbox { and } \frac {M}{c}\to \infty \mbox { as } c\to \infty .
\ee
Then
\be \label{lwrboundassy1}
F(c,M)>\exp \left \{M\left (1-\frac {\ln M+\ln \ln M}{c+1}+\frac {a^{\star }
+\ln c+\delta (c,M)}{c+1}
\right )\right \},
\ee
where
\be \label{deltacMlwr}
\delta (c,M)\to \infty \mbox { as } c\to \infty .
\ee
\end{Thm}

\mysection{An upper bound for Problem $Q_{c,M}$}\label{section7}
In this section we are concerned with the upper bound for $G(c,M)$ in (\ref{GcMsum}).
We will employ a method similar to the one used to establish a lower bound
for $F(c,M)$ in sections \ref{section5} and \ref{section6}.

It will be shown that the variables $G(c,M)$ satisfy relations similar to
(\ref{Frecurs}), and we wish to establish for $G(c,M)$ an inequality analogous to
(\ref{Fform}), with a reversed inequality sign. We note, however, that for fixed
$c$, $B$ and $\gamma $ the inequality
\be \label{Gform1}
G(c,M)\leq Be^{M\left (1-\frac {\ln M}{c+1}
+\frac {\gamma }{c+1}\right )}
\ee
cannot hold for sufficiently large $M$, since for such $M$ the right-hand side of
(\ref{Gform1}) becomes smaller than 1, while the left-hand side of (\ref{Gform1}) is
clearly larger than 1.

We henceforth focus on the function $G(c,M)$ defined in (\ref{GcMsum}).
Our goal is to estimate the value of $G(c,M)$ for pairs $(c,M)$ which belong to
the domain
$$
{\cal D}={\cal D}_{1/2}
$$
(recall (\ref{Dbetadfn})),
and we denote
\be \label{D1defntn}
{\cal D}_{+}=\{(c,M):e^{c/2}<M<e^{(c+1)/2}\}.
\ee
Analogous to
(\ref{Frecurs}), for points $(c,M)\in {\cal D}$ we have the following relation
\be \label{Grecurs}
G(c,M)=\sum _{z=0}^{[M/c]}G(c-1,M-cz)\frac {m_0^z}{z!}e^{z^2/m_0}.
\ee
(Of course, even though $(c,M)\in {\cal D}$, some points $(c-1,M-cz)$ in
(\ref{Grecurs}) may fail to belong to ${\cal D}$.)

To obtain an upper bound of the type (\ref{Gform1}) on ${\cal D}$  we will employ the
iterative method described in sections \ref{section5} and \ref{section6}. To use this
approach in the present situation we have to guarantee in advance that
(\ref{Gform1}) holds for points in ${\cal D}_{+}$. This property will follow from
Proposition \ref{upperbundln2} and the next result.
\begin{Prop}\label{G2Fcm}
The following relation holds:
\be \label{GFnuinq}
G(c,M)<2^cF(c,M).
\ee
\end{Prop}
{\em Proof}: We note that
$$
z\leq \frac {M}{c}\leq \frac {e^{(c+1)/2}}{c} \mbox { and } m_0>\frac {e^c}{c},
$$
implying
$$
\frac {z^2}{m_0}<\frac {e}{c}.
$$
It follows that $e^{z^2/m_0}<2$ if $c>e/\ln 2$. Now (\ref{GFnuinq}) follows from
(\ref{Frecurs}) and (\ref{Grecurs}), employing induction on $c$. $\hfill \Box$
\begin{Rem}\label{FreplG}
We will establish an upper bound for $F(c,M)$, and
then use $(\ref{GFnuinq})$ to estimate $G(c,M)$ from above. Thus we wish to establish
for $F$ an inequality of the form
\be \label{Fform2}
F(c,M)\leq Be^{M\left (1-\frac {\ln M}{c+1}
+\frac {\gamma }{c+1}\right )}
\ee
for some coefficient $B$ and a certain $\gamma $ (which may depend on $c$ and $M$),
and in view of $(\ref{GFnuinq})$ this will yield the estimate
\be \label{GdefnE}
G(c,M)\leq B\exp \left \{M\left (1-\frac {\ln M}{c+1}
+\frac {\gamma }{c+1}\right )+c\ln 2\right \}.
\ee
We note that under the assumption
\be \label{Mc2dvrg}
\frac {M}{c^2}\to \infty \mbox { as } c\to \infty,
\ee
the term $c\ln 2$ in the exponent in $(\ref{GdefnE})$ becomes negligible compared to the
other terms in the exponent when $c\to \infty $.
\end{Rem}

The following result is a consequence of Proposition \ref{upperbundln2}.
\begin{Prop}\label{insiD1}
Let ${\cal D}_{+}$ be as in $(\ref{D1defntn})$, and let $C$ be as in Proposition
$\ref{upperbundln2}$. Then $(\ref{Fform2})$, with $B=2$ and $\gamma =C$,
holds on ${\cal D}_{+}$.
\end{Prop}

We consider (\ref{Grecurs}) as a difference equation in ${\cal D}$ satisfying boundary
upper bounds on ${\cal D}_{+}$ as expressed in Proposition \ref{insiD1}.
For a fixed $\kappa >1$ let
$$
D_{\kappa }={\cal D}\cap \{1\leq c\leq \kappa \}
$$
which is a bounded set, and it follows that for any fixed $\gamma $,
$F(\cdot ,\cdot)$ satisfies (\ref{Fform2}) on $D_{\kappa }$
for some $B>1$ (depending on $\gamma $).

Suppose that we have an upper bound for $F(\cdot ,\cdot )$ on $D_{\kappa }$,
and we consider in the left hand side of (\ref{Frecurs}) pairs $(c,M)$ which belong to
$D_{\kappa +1}\setminus D_{\kappa }$. We will next show that for such $(c,M)$
the right hand side of (\ref{Frecurs}) involves pairs $(c-1,M-cz)$ for which an
upper bound of the form (\ref{Fform2}) has been already established. We will then
use these bounds to estimate the right hand side of (\ref{Frecurs}), thus establishing
an upper bound for $F(c,M)$.
\begin{Prop}\label{belong}
If $(c,M)\in D_{\kappa +1}\setminus D_{\kappa }$ then
\be \label{cMzblong}
(c-1,M-cz)\in D_{\kappa }\cup {\cal D}_{+}
\ee
for every $0\leq z\leq M/c$.
\end{Prop}
{\em Proof}: If $(c,M)\in D_{\kappa +1}$ then
$\displaystyle {M\leq e^{c/2}}$. Obviously
this can be written in the form
$$M\leq e^{\frac {(c-1)+1}{2}},$$
implying that $(c-1,M)\in {\cal D}_{+}$ if
$\displaystyle{M>e^{{(c-1)}/{2}}}$, and $(c-1,M)\in D_{\kappa }$
if $\displaystyle{M\leq e^{{(c-1)}/{2}}}$. $\hfill \Box$

It follows from Proposition \ref{belong} that each summand $F(c-1,M-cz)$
in the right hand side of (\ref{Grecurs}) may be bounded by employing a bound of the form
(\ref{Fform2}) for $(c-1,M-cz)$.

In analogy with (\ref{mzeE}) we have that
\be \label{GmzeE}
\frac {m_0^z}{z!}<e^{\bar E},
\ee
where similarly  to (\ref{Elogaritm})
\be \label{GElogaritm}
\bar E=\left (z\ln m_0-z\ln z+z\right ).
\ee
(In (\ref{GElogaritm}) we ignore the term $\sqrt z$ in (\ref{Stirapprx}), since we consider
now an upper bound)). Substituting $m_0=(e-1)e^c/c$ in (\ref{GElogaritm}) we obtain
$$
\bar E= cz-z\ln z+z(1+\ln (e-1))-z\ln c.
$$
Let $A$ be as in (\ref{Alogaritm}), and analogous to (\ref{c1Mzbnd}) we assume that
$$
F(c-1,M-cz)\leq Be^A,
$$
so that
$$
F(c-1,M-cz)\frac {m_0^z}{z!}\leq Be^{A+\bar E}.
$$
It follows that an upper bound for $A+\bar E$ is given by the function
$H(z)$ in (\ref{maxexpres}), where the variable $a$ (recall (\ref{adefinitn}))
is replaced by $a^{\star }$ in (\ref{limita}).
We still denote this function by $H(z)$, and
analogous to (\ref{FcMHD}) we have the relation
\be \label{GcMHD}
F(c-1,M-cz)\frac {m_0^z}{z!}<Be^{ H(z)}.
\ee

As in section \ref{section5}, we should maximize the function $ H(z)$ over
$0\leq z\leq [M/c]$. But in the present situation, since we are concerned
with an upper bound, we may use the maximum of $H(z)$ over the real
interval $0\leq z\leq M/c$ and do not have to restrict to the integers in this interval.

Summarizing the above discussion we obtain, analogous to (\ref{Hzmaxlower}), the
following result.
\begin{Prop}\label{GcMBbond}
Assume that
\be \label{GcMBinq}
F(c,M)\leq Be^{M\left (1-\frac {\ln M}{c+1}+\frac {\gamma}{c+1}\right )}
\ee
for every $(c,M)\in D_{\kappa }$, for some $\gamma>C$ and $\kappa >1$. Then
\be \label{maxestimate}
\max \left \{F(c-1,M-cz)\frac {m_0^z}{z!}:0\leq z\leq \frac {M}{c}\right \}
\leq Be^{M\left (1-\frac {\ln M}{c}+\frac {\gamma +f(\gamma )}{c}\right )},
\ee
implying
\be \label{muineqult}
F(c,M)<Be^{M\left (1-\frac {\ln M}{c}+\frac {\gamma +f(\gamma )}{c}\right )+\ln (M/c)}
\ee
for every $(c,M)\in D_{\kappa +1}$.
\end{Prop}
\begin{Rem}
The term $\ln (M/c)$ appears in $(\ref{muineqult})$ since we should multiply the
maximum in $(\ref{maxestimate})$ by the number of terms which appear in the sum in
$(\ref{Frecurs})$. We may use $\ln (M/c)$ rather than $\ln ([M/c]+1)$ since there are
in $(\ref{Frecurs})$ several summands which are much smaller than the maximal term there.
\end{Rem}
In this section we use induction to establish an inequality of the type
(\ref{Fform2}), with $\gamma $ depending on $(c,M)$ as follows:
\be \label{gamcM7}
\gamma (c,M)=\bar a+\ln c+\ln \ln c-\ln \ln M
\ee
for a certain $\bar a$.

{\bf The next is an important comment.}

We consider now the maximization
in the left hand side of (\ref{maxestimate}). Employing an induction hypothesis we obtain
bounds on the expressions $F(c-1,M-cz)$, using inequalities  of the form (\ref{GcMBinq})
for the pairs $(c-1,M')$, where $M'=M-cz$. In these bounds we denote
$\gamma =\gamma (c-1,M')$, using (\ref{gamcM7}).
Suppose that the {\em maximum over the bounds} is attained at $1<M_0\leq M$, and denote
$\gamma _0=\gamma (c-1,M_0)$, namely
\be \label{gama0defn}
\gamma _0=\bar a+\ln (c-1)+\ln \ln (c-1)-\ln \ln M_0.
\ee
Clearly the maximum over the bounds is not larger than the maximal value of
\be \label{gam0functn}
\exp \left \{(M-cz)\left [1-\frac {\ln (M-cz)}{c}+\frac {\gamma _0}{c}\right ]\right \}
\frac {m_0^z}{z!}
\ee
over $0\leq z\leq M/c$.

In view of (\ref{maxestimate}) and (\ref{muineqult}), and analogous to (\ref{Mgamampr1}),
we wish to establish
\be \label{gamafgamampr}
\frac {\gamma _0+f(\gamma _0)}{c}<\frac {\ln M}{c(c+1)}+\frac {\gamma '}{c+1},
\ee
where
\be \label{gamprdefn}
\gamma '=\bar a+\ln c+\ln \ln c-\ln \ln M.
\ee
We first address the term $f(\gamma _0)$ in (\ref{gamafgamampr}),
and recalling (\ref{fdefn}) we have
\be \label{fgamaapr}
f(\gamma _0)=\ln \left (1+e^{a^{\star }-\bar a}\frac {\ln M_0}{(c-1)\ln (c-1)}\right ).
\ee
We assume now that
$$
(c,M)\in {\cal D}_{\beta },
$$
and denote in (\ref{gamcM7})
\be \label{abardefn}
\bar a=a^{\star}+\delta
\ee
for some (not necessarily positive) $\delta $. It follows from (\ref{fgamaapr}) that
$$
f(\gamma _0)<e^{-\delta }\frac {\ln M}{(c-1)\ln (c-1)},
$$
and concerning (\ref{gamafgamampr})
we have thus established that
\be \label{fgamlnMest}
\frac {f(\gamma _0)}{c}<\frac {q\ln M}{c(c+1)\ln c}
\ee
for some constant $q>1$ independent of $\beta $ and $c$, if $c$ is sufficiently large.

We next consider the terms $\gamma _0/c$ and $\gamma '/(c+1)$ in (\ref{gamafgamampr}). Let
$z_0$ be the point where the maximization over $z$ of (\ref{gam0functn}) is attained, and
let, as above, $M_0=M-cz_0$. We note that in this maximization, the value $\gamma _0$ is
the same for all the points $(c-1,M')$, $1<M'\leq M$. We have then
\be \label{M0definit}
M_0=M(1-t_0),
\ee
where by (\ref{toexpr})
$$
t_0=\frac {1}{1+e^{\gamma _0-a^{\star }}}=\frac {q_1e^{-\delta }\ln M_0}{c\ln c},
$$
for some constant $q_1$ if $c$ is sufficiently large. Thus
\be \label{lnt0estm}
\ln (1-t_0)=-\frac {q_2\ln M_0}{c\ln c}
\ee
for some constant $q_2$, and it follows from (\ref{M0definit}) and
(\ref{lnt0estm}) that
$$
\left (1+\frac {q_2}{c\ln c}\right )\ln M_0=\ln M,
$$
hence
$$
\ln M_0=\left (1-\frac {q_3}{c\ln c}\right )\ln M
$$
for some $q_3>q_2$. The last relation implies that
\be \label{lnlnMM0}
\ln \ln M_0>\ln \ln M-\frac {2q_3}{c\ln c}
\ee
if $c$ is sufficiently large.

Using the expressions (\ref{gama0defn}) and  (\ref{gamprdefn}) it follows from
(\ref{lnlnMM0}) that
\newline
$\displaystyle{\frac {\gamma _0}{c}-\frac {\gamma '}{c+1}}$ is smaller than
$$
\frac {\bar a+\ln (c-1)+\ln \ln (c-1)}{c}
-\frac {\bar a+\ln c+\ln \ln c}{c+1}+\frac {2q_3}{c^2\ln c}
-\frac {\ln \ln M}{c(c+1)},
$$
implying that
\be \label{fracgam0gm2}
\frac {\gamma _0}{c}-\frac {\gamma '}{c+1}<
\frac {\bar a+\ln c+\ln \ln c}{c(c+1)}+\frac {2q_3}{c^2\ln c}
-\frac {\ln \ln M}{c(c+1)}.
\ee
Using $M>c$ we conclude from (\ref{fracgam0gm2}) that
\be \label{fracgam0gm3}
\frac {\gamma _0}{c}-\frac {\gamma '}{c+1}<\frac {\ln (kc)}{c(c+1)}
\ee
for large enough $c$, where we denote
$$
k=a^{\star }+1 .
$$

We next examine the inequality
\be \label{lnclnMinq}
\ln (kc)<\left [1-\frac {q}{\ln c}\right ]\ln M
\ee
where $q$ is as in (\ref{fgamlnMest}). We note that (\ref{gamafgamampr}) follows from
(\ref{fgamlnMest}), (\ref{fracgam0gm3}) and (\ref{lnclnMinq}), hence it only remains to
establish (\ref{lnclnMinq}). But (\ref{lnclnMinq}) holds if
\be \label{lnalphac}
\left (1+\frac {q_0}{\ln c}\right )\ln (kc)<\ln M
\ee
for a certain $q_0>q$, e.g. we may take $q_0=2q$ provided that $c$ satisfies $\ln c>2q$.
The inequality (\ref{lnalphac}), however, is equivalent to
$$
M>\left (e^{q_0}k^{1+q_0/\ln c}\right )c,
$$
which is satisfied if
\be \label{Mlmdac}
M>Kc
\ee
for the constant $K=e^{q_0}k^{1+q_0/\ln c}$. We have thus established the following result.
\begin{Prop}\label{discsummry}
Let the constant $K>1$ be fixed, and for some constant $\bar a$ let $\gamma (c,M)$
be as in $(\ref{gamcM7})$. Then there exist constants $B$ and $c_0$ such that
\be \label{acumineqlt}
F(c,M)<Be^{M\left (1-\frac {\ln M}{c+1}+\frac {\gamma (c,M)}{c+1}\right )+c\ln M}
\ee
holds provided that $c>c_0$.
\end{Prop}
{\em Proof}: The inequality (\ref{acumineqlt}) follows from (\ref{muineqult}) and
(\ref{gamafgamampr}) and the preceding discussion. We note that when employing successively
the inequalities (\ref{muineqult}) and (\ref{gamafgamampr}), the various terms
$\ln (M/c)$ in (\ref{muineqult}) accumulate, yielding the term $c\ln M$ in
(\ref{acumineqlt}). $\hfill \Box$

In the following result we consider pairs $(c,M)$ such that $M/c<K$.
\begin{Prop} \label{McsmalK}
There exist constants $B$ and $c_0$ such that the inequality
\be \label{Fuperform}
F(M,c)<Be^{M\left (1-\frac {\ln M}{c+1}+\frac {\gamma (c,M)}{c+1}\right )}
\ee
holds for every $(c,M)$ such that $1\leq M/c\leq K$ and $c>c_0$, where
$\gamma (c,M)$ is as in $(\ref{gamcM7})$.
\end{Prop}
{\em Proof}: We substitute $M=K_1c$ in (\ref{Fuperform}), for some $1\leq K_1\leq K$,
and use (\ref{gamcM7}) to obtain
\be \label{Fmlc}
F(c,M)<Be^{M\left (1-\frac {\ln K_1 -\bar a+2(\ln K_1)/(\ln c)}{c+1}\right )}.
\ee
But for $M=K_1c$ the right hand side of (\ref{Fmlc}) exceeds
$$
Be^Me^{-2K_1(\ln K_1-\bar a)}
$$
for large enough $c$. Since $F(c,M)<2N$ and $1\leq K_1 \leq K$, it follows that
(\ref{Fmlc}) indeed hold, provided that $B$ is sufficiently large. $\hfill \Box$

Propositions \ref{discsummry} and \ref{McsmalK} cover the whole range of interest, and
we summarize the above discussion as follows:
\begin{Prop}\label{upperbndn}
Let $\bar a$ be any constant. Then there exist constants $B$ and $c_0$ such that
\be \label{nineqult1}
F(c,M)<Be^{M\left (1-\frac {\ln M+\ln \ln M}{c+1}+\frac {\bar a+\ln c+\ln \ln c}
{c+1}\right )+c\ln M}
\ee
if $(c,M)\in {\cal D}$ and $c>c_0$.
\end{Prop}

We note that (\ref{nineqult1}) holds for arbitrarily small $\bar a$, for certain constants
$B$ and $c_0$ (depending on $\bar a$). This fact is due to the term
$\ln \ln c$ in the exponent in (\ref{nineqult1}).

Concerning $G(c,M)$, in view of Remark \ref{FreplG} we obtain the following results:
\begin{Prop}\label{upperbndG}
Let $\bar a$ be any constant. Then there exist constants $B$ and $c_0$ such that
\be \label{nineqult2}
G(c,M)<Be^{M\left (1-\frac {\ln M+\ln \ln M}{c+1}+\frac {\bar a+\ln c+\ln \ln c}
{c+1}\right )+c\ln 2M}
\ee
if $(c,M)\in {\cal D}$ and $c>c_0$.
\end{Prop}
\begin{Thm} \label{mainuprbndrslt}
Consider pairs $(c,M)$ such that
\be \label{uprcondi}
\frac {M}{c^2\ln M}\to \infty \mbox { as } c\to \infty .
\ee
Then there exists $c_0$ such that
\be \label{nineqult3}
G(c,M)<e^{M\left (1-\frac {\ln M+\ln \ln M}{c}+\frac {\bar a+\ln c+\ln \ln c}
{c}\right )}
\ee
if $(c,M)\in {\cal D}$  and $c>c_0$.
\end{Thm}

\mysection{The main results}\label{section8}
In this section we will establish our main results concerning lower and upper bounds
for $\nu (n,N)$. They consist of rephrasing the results in sections \ref{section6}
and \ref{section7} in terms of $n$ and $N$ instead of $c$ and $M$.

\begin{Prop}\label{noaymptot8}
Let $a^{\star}$ be defined by $(\ref{limita})$, and
let $a<a^{\star }$ be fixed. Then there exist $\beta >0$, $n_0$ and $b$ such that
\be \label{lwrboundassy8}
\frac {\ln \nu (n,N)}{\ln N}>1-\frac {\ln \ln N+\ln \ln \ln N}{\ln n}
+\frac {a+\ln \ln n}{\ln n}
\ee
for every $(n,N)$ such that $n^b<N<e^{n^{\beta }}$ and $n>n_0$.
\end{Prop}

Our first main result is concerned with the asymptotic lower bound for $\nu (n,N)$.

\begin{Thm} \label{mainlwrbndrslt8}
Consider pairs $(n,N)$ such that
\be \label{lwrcondi8}
\frac {\ln \ln N}{\ln n}\to 0 \mbox { and } \frac {\ln N}{\ln n}\to \infty
\mbox { as } n\to \infty .
\ee
Then
\be \label{lwrboundassy18}
\frac {\ln \nu (n,N)}{\ln N}>1-\frac {\ln \ln N+\ln \ln \ln N}
{\ln n}+\frac {a^{\star }
+\ln \ln n+\delta (n,N)}{\ln n},
\ee
where
\be \label{deltacMlwr8}
\delta (n,N)\to \infty \mbox { as } n\to \infty .
\ee
\end{Thm}

Concerning upper bounds for $\nu (n,N)$ we have the following result:

\begin{Prop}\label{upperbndG8}
Let $\bar a$ be any constant. Then there exist a constant $n_0$ such that
\be \label{nineqult28}
\frac {\ln \nu (n,N)}{\ln N}<1-\frac {\ln \ln N+\ln \ln \ln N}
{\ln n}+\frac {\bar a+\ln \ln n+\ln \ln \ln n}
{\ln n}+\frac {\ln n\ln (2\ln N)}{\ln N}
\ee
if $N<e^{\sqrt n}$ and $n>n_0$.
\end{Prop}

Our second main result is concerned with the asymptotic lower bound
for $\nu (n,N)$.

\begin{Thm} \label{mainuprbndrslt8}
Consider pairs $(n,N)$ such that
\be \label{uprcondi8}
\frac {\ln N}{(\ln n)^2\ln \ln N}\to \infty \mbox { as } n\to \infty .
\ee
Then there exists $n_0$ such that
\be \label{nineqult38}
\frac {\ln \nu (n,N)}{\ln N}<1-\frac {\ln \ln N+\ln \ln \ln N}{\ln n}
+\frac {\bar a+\ln \ln n+\ln \ln \ln n}{\ln n}
\ee
if $N<e^{\sqrt n}$ and $n>n_0$.
\end{Thm}

\mysection{Appendix}\label{section9}
{\em Proof of Theorem $\ref{niceresult}$}: Let $F=[1,N]\setminus E$ be the complement of $E$ in $[1,N]$. For a prime
$\sqrt N\leq p\leq N$ we denote by $F_p$ the set of integers in $F$ which are divisible
by $p$. Then $F_{p_1}\cap F_{p_2}=\emptyset $ if $p_1\not =p_2$,
$$
\#(F_p)=\left [\frac {N}{p}\right ]
$$
and it follows that
\be \label{Npsum}
\#(F)=\sum _{\sqrt N\leq p\leq N}[N/p]<N\sum _{p\geq \sqrt N}^N\frac {1}{p},
\ee
where the sum is over the primes in the indicated interval.
To estimate the sum in the right hand side of (\ref{Npsum}) we consider, more generally,
sums of the form
\be \label{anbsum}
S_{a,b}=\sum _{a\leq p\leq b}\frac {1}{p}.
\ee
By the Prime Numbers Theorem the distribution function of the number of primes in
the real line is, for large
enough $x$, $\Phi (x)=x/\ln x$.
Using this in the summation in (\ref{anbsum}) implies that for sufficiently large $a$
we have
$$
S_{a,b}\approx \int _a^b\frac {d\Phi (x)}{x}= \int _a^b\frac {\Phi(x)dx}{x^2}
+\left.\frac {\Phi (x)}{x}\right |_a^b,
$$
and substituting $\Phi (x)=x/\ln x$ we conclude that
\be \label{Sabsum}
S_{a,b}\approx \int _a^b\frac {dx}{x\ln x}
+\left.\frac {1}{\ln x}\right |_a^b <\ln \ln b-\ln \ln a.
\ee
For $a=\sqrt N$ and $b=N$ the right hand side of (\ref{Sabsum}) is equal to $\ln 2$,
and using this in (\ref{Npsum}) yields that for sufficiently large $N$ we have
$$\#(F)<N\ln 2,$$ implying $$\#(E)>N\ln (e/2).$$ This establishes (\ref{FsqrtN}) and
concludes the proof. $\hfill \Box$.

{\em Proof of Proposition} \ref{nlnNsqrbnd}:
It follows from $\nu (2,N)\leq \ln N/\ln 2$ that
$$
\nu (2,N)\leq \frac {\sqrt N}{\ln 2},
$$
since $\ln N<\sqrt N$ for every $N\geq 1$.
It is easy to see that
\be \label{nuprod}
\nu (p_k,N)\leq \frac {\sqrt N}{(\ln 2)(1-1/\sqrt {p_2})\cdots (1-1/\sqrt {p_k})},
\ee
for every $k\geq 2$. Relation (\ref{nuprod}) can be established by
employing a simple induction argument, using (\ref{nupkrecurs}).

To estimate from above the right hand side of (\ref{nuprod}), we have to estimate
from below the product
\be \label{prodestim}
\prod _{j=1}^k\left (1-\frac {1}{\sqrt {p_j}}\right ),
\ee
and for this we estimate from above the sum
\be \label{sumestim}
\sum _{j=1}^k\frac {1}{\sqrt {p_j}}.
\ee
To this end we use the distribution function
$$
\Phi (x)=\frac {x}{\ln x}
$$
of the primes in the real line, and we have to estimate
$$
\int _3^{p_k }\frac {d\Phi (x)}{\sqrt x}.
$$
This leads to
\be \label{CpklnlnN}
\int _3^{p_k }\frac {dx}{\sqrt x\ln x}=\int _{\sqrt 3}^{\sqrt {p_k}}\frac {dt}{2\ln t}<\frac {C\sqrt {p_k}}{\ln p_k}
\ee
for some constant $C>0$, and we obtain
\be \label{nupkNsq}
\nu (p_k,N)\leq N^{1/2}e^{\frac {C\sqrt {p_k}}{\ln p_k}}.
\ee
For a prescribed $n=\alpha (\ln N)^2$
we let $p_k$ be the smallest prime $p$ which satisfies $p\geq n$. Employing
(\ref{nupkNsq}) for this $p_k$ yields
the assertion of the proposition.
$\hfill \Box$
\end{document}